\documentstyle[12pt]{article}
\begin{document}

\title{Some Inequalities satisfied by Periodical \\
Solutions of Multi-Time Hamilton Equations}
\author{Iulian Duca, Constantin Udri\c ste}
\date{}
\maketitle

\begin{abstract}
The objective of this paper is to find some inequalities satisfied by
periodical solutions of multi-time Hamilton systems, when the Hamiltonian is
convex. To our knowledge, this subject of first-order field theory is still
open.

Section 1 recall well-known facts regarding the equivalence between
Euler-Lagrange equations and Hamilton equations and analyses the action that
produces multi-time Hamilton equations, emphasizing the role of the
polysymplectic structure. Section 2 extends two inequalities of [21] from a
cube to parallelipiped and proves two inequalityes concerning multiple
periodical solutions of multi-time Hamilton equations.
\end{abstract}



\section{Multi-time Hamilton equations and polysymplectic structure}

The paper studies the solutions with multiple periodicity of the Hamilton
multi-time equations.

A function $u=\left( u^{1},...,u^{n}\right) $ with many variables $\left(
t^{1},...,t^{p}\right) $, is multiple periodical with the period $T=\left(
T^{1},...,T^{p}\right) \in R^{p}$ if 
\[
u\left( t^{1}+k_{1}T^{1},...,t^{p}+k_{p}T^{p}\right) =u\left(
t^{1},...,t^{p}\right) , 
\]
where $k_{1},...,k_{p}$ are integers. We consider the function $u$ defined
on the parallelepiped $T_{0}=\left[ 0,T^{1}\right] \times \left[ 0,T^{2}%
\right] \times ...\times \left[ 0,T^{p}\right] \subset R^{p}$, with values
in $R^{n}$. We will denote by $T=\left( T^{1},...,T^{p}\right) \in R^{p}$.\
The existence of the weak gradient of the function $u$ assures the multiple
periodicity of the function $u$. We use the Hilbert space $H_{T}^{1}$
attached to the Sobolev space $W_{T}^{1,2}$ of the functions $u\in
L^{2}\left( T_{0},R^{n}\right) $ which have a weak gradient $\displaystyle%
\frac{\partial u}{\partial t}\in L^{2}\left( T_{0},R^{n}\right) $. The
Wirtinger inequality from this paper has a specific form because of the
multidimensional character of the definition domain $T_{0}$. The
inequalities from theorems 3 and 4 constitute generalizations of some
theorems of [5], from the particular case $p=1$ to an arbitrary $p$.

The Euclidean structure on $R^n$ is based on the scalar product $\left(
u,v\right) =\delta _{ij}u^{i}v^{j}$, and the norm $\left| u\right| =\sqrt{%
\delta _{ij}u^{i}u^{j}}$. The Hilbert space $H_{T}^{1}$ is endowed with the
scalar product

\[
\left\langle u,v\right\rangle =\int_{T_{0}}\left( \delta _{ij}u^{i}\left(
t\right) v^{j}\left( t\right) +\delta _{ij}\delta ^{\alpha \beta }\frac{%
\partial u^{i}}{\partial t^{\alpha }}\left( t\right) \frac{\partial v^{j}}{%
\partial t^{\beta }}\left( t\right) \right) dt^{1}\wedge ...\wedge dt^{p}, 
\]
and the corresponding norm $\sqrt{\left\langle u,u\right\rangle }=\left\|
u\right\| .$

\subsection{Multi-time Hamilton equations}

We consider the multi-time variable $t=\left( t^{1},...,t^{p}\right) \in
T_{0}\subset R^{p}$, the functions $x^{i}:R^{p}\rightarrow R$ , $\left(
t^{1},...,t^{p}\right) \rightarrow x^{i}\left( t^{1},...,t^{p}\right) $, $%
i=1,...,n$ , and the partial velocities $x_{\alpha }^{i}=\displaystyle\frac{%
\partial x^{i}}{\partial t^{\alpha }}$ , $\alpha =1,...,p$.

{\bf Definition 1 }The PDEs 
\[
\frac{\partial }{\partial t^{\alpha }}\frac{\partial L}{\partial x_{\alpha
}^{i}}=\frac{\partial L}{\partial x^{i}},\quad i=1,...,n,\quad \alpha
=1,...,p 
\]
(second order PDEs system on the n-dimensional space) are called
Euler-Lagrange equations for the Lagrangian 
\[
L:R^{p+n+np}\rightarrow R,\quad \left( t^{\alpha },x^{i},x_{\alpha
}^{i}\right) \rightarrow L\left( t^{\alpha },x^{i},x_{\alpha }^{i}\right) 
\]

The Hamilton equations in the multi-time case are obtained using the partial
derivatives (polymomenta) 
$$
p_{k}^{\alpha }=\frac{\partial L}{\partial x_{\alpha }^{k}}\eqno(1) 
$$
and the Hamiltonian $H=p_{k}^{\alpha }x_{\alpha }^{k}-L$. If \ $L$ satisfies
some regularity conditions, then the system (1) defines a $C^{1}$ bijective
transformation $x_{\alpha }^{i}\rightarrow p_{i}^{\alpha }$ , called the
Legendre transformation for the multi-time case. By this transformation we
have 
\[
\frac{\partial H}{\partial p_{i}^{\alpha }}=x_{\alpha }^{i}+p_{k}^{\beta }%
\frac{\partial x_{\beta }^{k}}{\partial p_{i}^{\alpha }}-\frac{\partial L}{%
\partial x_{\beta }^{k}}\frac{\partial x_{\beta }^{k}}{\partial
p_{i}^{\alpha }}=x_{\alpha }^{i} 
\]
\[
\frac{\partial H}{\partial x^{i}}=p_{k}^{\alpha }\frac{\partial x_{\alpha
}^{k}}{\partial x^{i}}-\frac{\partial L}{\partial x^{i}}-\frac{\partial L}{%
\partial x_{\alpha }^{k}}\frac{\partial x_{\alpha }^{k}}{\partial x^{i}}=-%
\frac{\partial L}{\partial x^{i}}. 
\]
Consequently, the $np+n$ Hamilton equations 
\[
\displaystyle\frac{\partial x^{i}}{\partial t^{\alpha }}=\displaystyle\frac{%
\partial H}{\partial p_{i}^{\alpha }}, 
\]
\[
\displaystyle\frac{\partial p_{i}^{\alpha }}{\partial t^{\alpha }}=-%
\displaystyle\frac{\partial H}{\partial x^{i}} 
\]
(summation after $\alpha $), $i=1,...,n$, $\alpha =1,...,p$ are first order
PDEs on the space $R^{n+pn}$, equivalent to the Euler-Lagrange equations on $%
R^{n}$.

There are different point of views to study these equations which appear in
first-order field theory (see [1]-[3], [8]-[10], [12]-[22]). In our context,
we need of Hilbert-Sobolev space methods for PDEs ([4], [6], [11]).

Let us write the multi-time Hamilton equations in the form 
\[
\displaystyle\delta _{\beta }^{\alpha }\delta _{j}^{i}\frac{\partial
p_{i}^{\beta }}{\partial t^{\alpha }}+\frac{\partial H}{\partial x^{j}}%
=0,\quad \displaystyle-\delta _{\beta }^{\alpha }\delta _{j}^{i}\frac{%
\partial x^{j}}{\partial t^{\alpha }}+\frac{\partial H}{\partial
p_{i}^{\beta }}=0, i, j=1,...,n; \alpha, \beta =1,...,p 
\]
or 
$$
\left( \delta \otimes J\right)\displaystyle\frac{\partial u}{\partial t} = -
\nabla H, \eqno(2) 
$$
where 
\[
\delta \otimes J=\left( 
\begin{array}{cc}
0 & \delta _{\beta }^{\alpha }\delta _{j}^{i} \\ 
\noalign{\medskip}-\delta _{\beta }^{\alpha }\delta _{j}^{i} & 0
\end{array}
\right), \displaystyle\frac{\partial u}{\partial t}=\left( 
\begin{array}{c}
\displaystyle\frac{\partial x^{j}}{\partial t^{\alpha }} \\ 
\noalign{\medskip}\displaystyle\frac{\partial p_{i}^{\beta }}{\partial
t^{\alpha }}
\end{array}
\right) , \nabla H =\left( 
\begin{array}{c}
\displaystyle\frac{\partial H}{\partial x^{j}} \\ 
\noalign{\medskip}\displaystyle\frac{\partial H}{\partial p_{j}^{\beta }}
\end{array}
\right). 
\]
The operator $\delta \otimes J$ is a polysymplectic structure acting on $%
R^{np+np^{2}}$ with values in \ $R^{n+np}.$ The (1,2)-block $\delta _{\beta
}^{\alpha }\delta _{j}^{i} $ acts linearly by $\delta _{j}^{i}$ and tracely
by $\delta _{\beta }^{\alpha }$. The (2,1)-block $-\delta _{\beta }^{\alpha
}\delta _{j}^{i} $ acts lineraly both by $\delta _{j}^{i}$ and $\delta
_{\beta }^{\alpha }$. The operator $\delta \otimes J$ induces a
multisymplectic PDE operator $(\delta \otimes J)\displaystyle\frac{\partial}{%
\partial t}$ which work as follows 
\[
(\delta \otimes J)\displaystyle\frac{\partial}{\partial t}\left( 
\begin{array}{c}
x \\ 
p
\end{array}
\right): \left( 
\begin{array}{cc}
0 & \delta _{\beta }^{\alpha }\delta _{j}^{i} \\ 
\noalign{\medskip}-\delta _{\beta }^{\alpha }\delta _{j}^{i} & 0
\end{array}
\right)\displaystyle\frac{\partial}{\partial t^\alpha}\left( 
\begin{array}{c}
x^j \\ 
p^{\beta}_i
\end{array}
\right)=\left( 
\begin{array}{c}
\displaystyle\frac{\partial p^{\alpha}_{j}}{\partial t^{\alpha }} \\ 
\noalign{\medskip}\displaystyle -\frac{\partial x^{j}}{\partial t^{\beta }}
\end{array}
\right). 
\]
Repeating we obtain the square

\[
(\delta \otimes J)\displaystyle\frac{\partial}{\partial t}\left( 
\begin{array}{c}
\hbox {div} p \\ 
-\displaystyle \frac{\partial x}{\partial t}
\end{array}
\right): \left( 
\begin{array}{cc}
0 & \delta _{\beta }^{\alpha }\delta _{j}^{i} \\ 
\noalign{\medskip}-\delta _{\beta }^{\alpha }\delta _{j}^{i} & 0
\end{array}
\right)\left( 
\begin{array}{c}
\displaystyle\frac{\partial^2 p^{\gamma}_i}{\partial t^{\alpha}\partial
t^{\gamma}} \\ 
-\displaystyle\frac{\partial^2 x^j}{\partial t^{\alpha}\partial t{\beta}}
\end{array}
\right)=\left( 
\begin{array}{c}
-\Delta x^i \\ 
-\displaystyle\frac{\partial^2 p^{\gamma}_i}{\partial t^{\beta}\partial
t^{\gamma}}
\end{array}
\right). 
\]

\subsection{The action that produces \ multi-time Hamilton equations}

We consider a Hamiltonian $H:T_{0}\times R^{n}\times R^{np}\rightarrow
R,\left( t,u\right) \rightarrow H\left( t,u\right) $ whose restriction $%
H\left( t,\cdot \right) $ is $C^{1}$ and convex.

{\bf Theorem 1 }[21]{\bf \ }{\it Let $u=\left( x,p\right)$. The action \ }$%
{\it \Psi }${\it , whose Euler-Lagrange equations are the Hamilton
equations, is } 
\[
{\it \Psi \left( u\right) =\int_{T_{0}}}{\cal L}{\it \left( t,u,\frac{%
\partial u}{\partial t}\right) dt^{1}\wedge ...\wedge dt^{p},\quad } 
\]

\[
{\cal L}\left( t,u,\frac{\partial u}{\partial t}\right) = =-\frac{1}{2}%
G\left( \delta \otimes J\frac{\partial u}{\partial t},u\right) -H\left(
t,u\right) , 
\]

{\it where the scalar product is represented by the matrix } 
\[
{\it G=\left( 
\begin{array}{cc}
\delta _{ij} & 0 \\ 
\noalign{\medskip}0 & \delta ^{\beta \alpha }\delta _{ij}
\end{array}
\right) } 
\]
{\it \ (standard Riemannian metric from }${\it R^{n+np}}${\it ). }

{\bf Proof}. Indeed, the Euler-Lagrange equations produced by 
\[
{\cal L}=-\frac{1}{2}\left( \frac{\partial p_{i}^{\alpha }}{\partial
t^{\alpha }}x^{i}-\frac{\partial x^{i}}{\partial t^{\alpha }}p_{i}^{\alpha
}\right) -H\left( t,x,p\right) = 
\]

\[
=-\frac{1}{2}\left( \frac{\partial p_{i}^{\alpha }}{\partial t^{\alpha }},-%
\frac{\partial x^{j}}{\partial t^{\beta }}\right) \left( 
\begin{array}{cc}
\delta _{ij} & o \\ 
\noalign{\medskip}o & \delta ^{\alpha \beta }\delta _{ij}
\end{array}
\right) \left( 
\begin{array}{c}
x^{i} \\ 
\noalign{\medskip}p_{i}^{\alpha }
\end{array}
\right) -H\left( t,x\left( t\right) ,p\left( t\right) \right) 
\]

can be rewritten 
\[
\frac{1}{2}\frac{\partial p_{i}^{\alpha }}{\partial t^{\alpha }}=-\frac{1}{2}%
\frac{\partial p_{i}^{\alpha }}{\partial t^{\alpha }}-\frac{\partial H}{%
\partial x^{i}}\;\hbox{,\,\,i.e.\,\,,}\;\frac{\partial p_{i}^{\alpha }}{%
\partial t^{\alpha }}=-\frac{\partial H}{\partial x^{i}} 
\]
and 
\[
\frac{\partial x^{i}}{\partial t^{\alpha }}=\frac{\partial H}{\partial
p_{i}^{\alpha }}. 
\]

\section{Basic inequalities}

\subsection{\bf \ Wirtinger multi-time inequality}

In $L^{2}\left( T_{0},R^{n}\right) $ we use the scalar product 
\[
\langle u,v\rangle =\int_{T_{0}}\left( \delta _{ij}u^{i}v^{j}\right)
dt^{1}\wedge ...\wedge dt^{p} 
\]
and the norm $\left\| u\right\| _{L^{2}}=\sqrt{\langle u,u\rangle }.$
Similarly, in $L^{2}\left( T_{0},C^{n}\right) $ we use the scalar product 
\[
\langle u,v\rangle =\int_{T_{0}}\left( \delta _{ij}u^{i}{\overline v}%
^{j}\right) dt^{1}\wedge ...\wedge dt^{p} 
\]
and the norm $\left\| u\right\| _{L^{2}}=\sqrt{\langle u,{\overline u}%
\rangle }.$

Let us extend the Theorem 4.4 from [21] to the parallelipiped $T_{0}.$

{\bf Theorem 2 }\bigskip {\it Any function }${\it u}${\it \ from }${\it H}_{%
{\it T}}^{{\it 1}}${\it \ with mean zero satisfies the inequality} 
\[
{\it \displaystyle\int_{T_{0}}\left| u\left( t\right) \right|
^{2}dt^{1}\wedge ...\wedge dt^{p}\leq \frac{\left( \displaystyle%
\max_{i}\left\{ T^{i}\right\} \right) ^{2}}{4\pi ^{2}}\int_{T_{0}}\left| 
\frac{\partial u}{\partial t}\right| ^{2}dt^{1}\wedge ...\wedge dt^{p}.}
\]
{\it \ }

{\bf Proof}. We express the function $u$ as the sum of a multiple Fourier
series 
\[
u\left( t\right) =\left( u^{1}\left( t^{1},...,t^{p}\right) ,...,u^{n}\left(
t^{1},...,t^{p}\right) \right) 
\]

\begin{eqnarray*}
&=&\left( \displaystyle\sum_{j_{1}^{1},...,j_{1}^{p}\in
Z^{p}}^{{}}C_{j_{1}^{1},...,j_{1}^{p}}^{1}e^{i\left( \frac{2\pi }{T^{1}}%
\cdot j_{1}^{1}t^{1}+...+\frac{2\pi }{T^{p}}\cdot j_{1}^{p}t^{p}\right)
},\right. \\
&&...,\left. \displaystyle\sum_{j_{n}^{1},...,j_{n}^{p}\in
Z^{p}}^{{}}C_{j_{n}^{1},...,j_{n}^{p}}^{n}e^{i\left( \frac{2\pi }{T^{1}}%
\cdot j_{n}^{1}t^{1}+...+\frac{2\pi }{T^{p}}\cdot j_{n}^{p}t^{p}\right)
}\right) .
\end{eqnarray*}
We calculate the square of the norm 
\[
\left( u,u\right) =\int_{T_{0}}\left( u\left( t\right) ,\overline{u\left(
t\right) }\right) dt^{1}\wedge ...\wedge dt^{p}= 
\]

\[
=\displaystyle\sum_{j_{1}^{1},...,j_{1}^{p}\in Z^{p}}^{{}}\left(
C_{j_{1}^{1},...,j_{1}^{p}}^{1}\right) ^{2}\int_{0}^{T^{1}}e^{i\frac{2\pi }{%
T^{1}}\left( j_{1}^{1}-j_{1}^{1}\right) t^{1}}dt^{1}...\int_{0}^{T^{p}}e^{i%
\frac{2\pi }{T^{p}}\left( j_{1}^{p}-j_{1}^{p1}\right) t^{p}}dt^{p}+... 
\]

\[
...+\displaystyle\sum_{j_{n}^{1},...,j_{n}^{p}\in Z^{p}}^{{}}\left(
C_{j_{n}^{1},...,j_{n}^{p}}^{1}\right) ^{2}\int_{0}^{T^{1}}e^{i\frac{2\pi }{%
T^{1}}\left( j_{n}^{1}-j_{n}^{1}\right) t^{1}}dt^{1}...\int_{0}^{T^{p}}e^{i%
\frac{2\pi }{T^{p}}\left( j_{n}^{p}-j_{n}^{p}\right) t^{p}}dt^{p} 
\]

\[
=\displaystyle\sum_{j_{1}^{1},...,j_{1}^{p}\in Z^{p}}^{{}}\left(
C_{j_{1}^{1},...,j_{1}^{p}}^{1}\right) ^{2}T^{1}...T^{p}+...+\displaystyle%
\sum_{j_{1}^{1},...,j_{1}^{p}\in Z^{p}}^{{}}\left(
C_{j_{n}^{1},...,j_{n}^{p}}^{n}\right) ^{2}T^{1}...T^{p}. 
\]
If we denote 
\[
C_{k_{1},...,k_{p}}=\left(
C_{k_{1},...,k_{p}}^{1},...,C_{k_{1},...,k_{p}}^{n}\right) 
\]
and 
\[
u=\displaystyle\sum_{\left( k_{1},...,k_{p}\right) \in
Z^{p}}^{{}}C_{k_{1},...,k_{p}}e^{i2\pi \left( \frac{k_{1}}{T^{1}}t_{1}+...+%
\frac{k_{p}}{T^{p}}t_{p}\right) }, 
\]
we find 
\[
\left( u,u\right) =\displaystyle\sum_{\left( k_{1},...,k_{p}\right) \in
Z^{p}}^{{}}\left| C_{k_{1},...,k_{p}}\right| ^{2}T^{1}...T^{p}. 
\]
Similarly, we consider the scalar product

\[
\left( \frac{\partial u}{\partial t},\frac{\partial v}{\partial t}\right)
=\int_{T_{0}}\delta _{ij}\delta ^{\alpha \beta }\frac{\partial u^{i}}{%
\partial t^{\alpha }}\overline{\frac{\partial v^{j}}{\partial t^{\beta }}}%
dt^{1}\wedge ...\wedge dt^{p}. 
\]
It follows the square of the norm 
\[
\left( \frac{\partial u}{\partial t},\frac{\partial u}{\partial t}\right) 
\]

\[
=\displaystyle\sum_{\left( k_{1},...,k_{p}\right) \in
Z^{p}}^{{}}T^{1}...T^{p}\left[ \left( C_{k_{1}...k_{p}}^{1}\right)
^{2}\left( \frac{2\pi k_{1}}{T^{1}}\right) ^{2}+...+\left(
C_{k_{1}...k_{p}}^{1}\right) \left( \frac{2\pi k_{p}}{T^{p}}\right)
^{2}+...\right. 
\]

\[
\left. ...+\left( C_{k_{1}...k_{p}}^{n}\right) ^{2}\left( \frac{2\pi k_{1}}{%
T^{1}}\right) ^{2}+...+\left( C_{k_{1}...k_{p}}^{n}\right) ^{2}\left( \frac{%
2\pi k_{p}}{T^{p}}\right) ^{2}\right] 
\]

\[
=\displaystyle\sum_{\left( k_{1},...,k_{p}\right) \in
Z^{p}}^{{}}T^{1}...T^{p}\left[ \left| C_{k_{1}...k_{p}}\right| ^{2}\left( 
\frac{2\pi k_{1}}{T^{1}}\right) ^{2}+...+\left| C_{k_{1}...k_{p}}\right|
^{2}\left( \frac{2\pi k_{p}}{T^{p}}\right) ^{2}\right] 
\]

\[
=\displaystyle\sum_{\left( k_{1},...,k_{p}\right) \in
Z^{p}}^{{}}T^{1}...T^{p}\left[ \left| C_{k_{1}...k_{p}}\right| ^{2}4\pi
^{2}\left( \left( \frac{k_{1}}{T^{1}}\right) ^{2}+...+\left( \frac{k_{p}}{%
T^{p}}\right) ^{2}\right) \right] 
\]

\[
\geq \frac{4\pi ^{2}}{\left( \displaystyle\max_{i}\left\{ T^{i}\right\}
\right) ^{2}}\displaystyle\sum_{\left( k_{1},...,k_{p}\right) \in
Z^{p}}^{{}}T^{1}...T^{p}\left| C_{k_{1}...k_{p}}\right| ^{2}\left(
k_{1}^{2}+...+k_{1}^{2}\right) 
\]

\[
\geq \frac{4\pi ^{2}}{\left( \displaystyle\max_{i}\left\{ T^{i}\right\}
\right) ^{2}}\displaystyle\sum_{\left( k_{1},...,k_{p}\right) \in
Z^{p}}^{{}}T^{1}...T^{p}\left| C_{k_{1}...k_{p}}\right| ^{2} 
\]
\[
\geq \frac{4\pi ^{2}}{\left( \displaystyle\max_{i}\left\{ T^{i}\right\}
\right) ^{2}}\int_{T_{0}}\left| u\left( t\right) \right| ^{2}dt^{1}\wedge
...\wedge dt^{p}. 
\]
Consequently 
\[
\int_{T_{0}}\left| u\left( t\right) \right| ^{2}dt^{1}\wedge ...\wedge
dt^{p}\leq \frac{\left( \displaystyle\max_{i}\left\{ T^{i}\right\} \right)
^{2}}{4\pi ^{2}}\int_{T_{0}}\left| \frac{\partial u}{\partial t}\right|
^{2}dt^{1}\wedge ...\wedge dt^{p}, 
\]
and this ends the proof.

\subsection{An estimate of the quadratic form \newline
$\displaystyle\int_{T_{0}}\left( \protect\delta \otimes J\frac{\partial u}{%
\partial t},u\left( t\right) \right) dt^{1}\wedge ...\wedge dt^{p}$}

Let us extend the Theorem 4.5 from [21] to the parallelipiped $T_{0}.$

{\bf Theorem 3 }{\it For any }${\it u\in H}_{{\it T}}^{{\it 1}}${\it \ we
have}

\[
{\it \displaystyle\int_{T_{0}}\left( \delta \otimes J\frac{\partial u}{%
\partial t},u\left( t\right) \right) dt^{1}\wedge ...\wedge dt^{p}} 
\]
\[
{\it \geq -\frac{\sqrt{p}\displaystyle\max_{i}\left\{ T^{i}\right\} }{2\pi }%
\int_{T_{0}}\left| \frac{\partial u}{\partial t}\right| ^{2}dt^{1}\wedge
...\wedge dt^{p}.} 
\]

{\bf Proof}. We denote $\widetilde{u}\left( t\right) =u\left( t\right) -%
\displaystyle\int_{T_{0}}u\left( t\right) dt^{1}\wedge ...\wedge dt^{p}$..
By using the Cauchy-Schwarz inequality and the multiple periodicity of $u$
we obtain the inequality

\[
\int_{T_{0}}\left( \delta \otimes J\frac{\partial u}{\partial t},u\left(
t\right) \right) dt^{1}\wedge ...\wedge dt^{p} 
\]
\[
=\int_{T_{0}}\left( \delta \otimes J\frac{\partial u}{\partial t},\widetilde{%
u}\left( t\right) +\int_{T_{0}}u\left( t\right) dt^{1}\wedge ...\wedge
dt^{p}\right) dt^{1}\wedge ...\wedge dt^{p} 
\]

\[
=\int_{T_{0}}\left( \delta \otimes J\frac{\partial u}{\partial t},\widetilde{%
u}\left( t\right) \right) dt^{1}\wedge ...\wedge dt^{p} 
\]
\[
+\int_{T_{0}}\left( \delta \otimes J\frac{\partial u}{\partial t}%
,\int_{T_{0}}u\left( t\right) dt^{1}\wedge ...\wedge dt^{p}\right)
dt^{1}\wedge ...\wedge dt^{p} 
\]

\[
=\int_{T_{0}}\left( \delta \otimes J\frac{\partial u}{\partial t},\widetilde{%
u}\left( t\right) \right) dt^{1}\wedge ...\wedge dt^{p} 
\]

\[
\geq -\left( \int_{T_{0}}\left| \delta \otimes J\frac{\partial u}{\partial t}%
\right| ^{2}dt^{1}\wedge ...\wedge dt^{p}\right) ^{\frac{1}{2}}\left(
\int_{T_{0}}\left| \tilde{u}\left( t\right) \right| ^{2}dt^{1}\wedge
...\wedge dt^{p}\right) ^{\frac{1}{2}} 
\]
From the inequality given by the Theorem 2 we have

\[
\int_{T_{0}}\left| \widetilde{u}\left( t\right) \right| ^{2}dt^{1}\wedge
...\wedge dt^{p}\leq \frac{\left( \displaystyle\max_{i}\left\{ T^{i}\right\}
\right) ^{2}}{4\pi ^{2}}\int_{T_{0}}\left| \frac{\partial \widetilde{u}}{%
\partial t}\right| ^{2}dt^{1}\wedge ...\wedge dt^{p} 
\]
\[
=\frac{\left( \displaystyle\max_{i}\left\{ T^{i}\right\} \right) ^{2}}{4\pi
^{2}}\int_{T_{0}}\left| \frac{\partial u}{\partial t}\right|
^{2}dt^{1}\wedge ...\wedge dt^{p}. 
\]
Because $\left| \delta \otimes J\displaystyle\frac{\partial u}{\partial t}%
\right| ^{2}\leq p\left| \displaystyle\frac{\partial u}{\partial t}\right|
^{2}$, we obtain 
\[
\int_{T_{0}}\left( \delta \otimes J\frac{\partial u}{\partial t},u\left(
t\right) \right) dt^{1}\wedge ...\wedge dt^{p} 
\]

\[
\geq -\left( \int_{T_{0}}\left| \delta \otimes J\frac{\partial u}{\partial t}%
\right| ^{2}dt^{1}\wedge ...\wedge dt^{p}\right) ^{\frac{1}{2}}\frac{%
\displaystyle\max_{i}\left\{ T^{i}\right\} }{2\pi } 
\]

\[
\cdot \left( \int_{T_{0}}\left| \frac{\partial u}{\partial t}\right|
^{2}dt^{1}\wedge ...\wedge dt^{p}\right) ^{\frac{1}{2}} 
\]
\[
\geq -\sqrt{p}\,\,\, \frac{\displaystyle\max_{i}\left\{ T^{i}\right\} }{2\pi 
}\int_{T_{0}}\left| \frac{\partial u}{\partial t}\right| ^{2}dt^{1}\wedge
...\wedge dt^{p}. 
\]

\subsection{Inequalities satisfied by periodical solutions of multi-time
Hamilton equations}

Let us find properties of solutions of $\left( \delta \otimes J\right) %
\displaystyle\frac{\partial u}{\partial t}+\nabla H\left( t,u\left( t\right)
\right) =0$ a.e. on $T_{0}$ satisfying the boundary conditions 
\[
u\mid _{S^{+}}=u\mid _{S^{-}}, 
\]
were $S^{+}$\ and $S^{-}$\ are opposite sides of the parallelipiped $T_{0}.$
Practically, we refer to bounds for such solutions.

{\bf Theorem 4 }{\it We consider the Hamiltonian\ } 
\[
{\it H:T}_{{\it 0}}{\it \times R}^{{\it n+np}}{\it \rightarrow R,\left(
t,u\right) \rightarrow H\left( t,u\right) } 
\]
{\it like a measurable function in }${\it t}${\it \ for any }${\it u\in R}^{%
{\it n+np}}$, and $C^{1}$ {\it convex in }${\it u}${\it \ for any} 
\[
{\it \ t\in T_{0}=\left[ 0,T^{1}\right] \times ...\times \left[ 0,T^{p}%
\right] \subset R}^{{\it p}}{\it .} 
\]
{\it If there exists} {\it the constants} 
\[
{\it \alpha \in \left( 0,\frac{\pi }{\sqrt{p}\displaystyle\max_{i}\left\{
T^{i}\right\} }\right) ,\beta \geq 0,\gamma \geq 0,\delta \geq 0} 
\]
{\it such that } 
\[
{\it \delta \left| u\right| +\beta \leq H\left( t,u\right) \leq \frac{\alpha 
}{2}\left| u\right| ^{2}+\gamma } 
\]
{\it for all }${\it t\in T}_{{\it 0}}${\it \ and }${\it u\in R}^{{\it n+np}}$%
,{\it \ then, any multiple periodical solution } 
\[
{\it u=\left( x_{i},p_{i}^{\alpha }\right) ,i=1,...,n,\alpha =1,...,p,} 
\]
{\it of the equation} 
$$
{\it \delta \otimes J\frac{\partial u}{\partial t}+\nabla H\left( t,u\left(
t\right) \right) =0,}\eqno(3) 
$$
verifies the inequalities {\it \ } 
$$
{\it \displaystyle\int_{T_{0}}\left| \frac{\partial u}{\partial t}\right|
^{2}dt^{1}\wedge ...\wedge dt^{p}\leq \frac{2\alpha \left( \beta +\gamma
\right) \pi T^{1}...T^{p}}{\pi -\alpha \displaystyle\max_{i}\left\{
T^{i}\right\} \sqrt{p}}}\eqno(4) 
$$
$$
{\it \displaystyle\int_{T_{0}}\left| u\left( t\right) \right| dt^{1}\wedge
...\wedge dt^{p}\leq \frac{\pi T^{1}...T^{p}\left( \beta +\gamma \right) }{%
\delta \left( \pi -\alpha \displaystyle\max_{i}\left\{ T^{i}\right\} \sqrt{p}%
\right) }.}\eqno(5) 
$$
{\bf Proof. } From the inequality 
\[
\delta \left| u\right| ^{2}-\beta \leq H\left( t,u\right) \leq \displaystyle%
\frac{\alpha }{2}\left| u\right| ^{2}+\gamma 
\]
we obtain 
\[
-\beta \leq H\left( t,u\right) \leq \alpha 2^{-1}\left| u\right| ^{2}+\gamma
. 
\]
By applying [5, Proposition 2.2], considering $F\left( u\right) =H\left(
t,u\right) $, $p=q=2$, $v=\nabla H\left( t,u\right) $ we obtain 
\[
\displaystyle\frac{1}{2\alpha }\left| \nabla H\left( t,u\right) \right|
^{2}\leq \left( \nabla H\left( t,u\right) ,u\right) +\beta +\gamma . 
\]
Because $u$ is the solution of the equation (3), we have $\nabla H\left(
t,u\right) =-\delta \otimes J\displaystyle\frac{\partial u}{\partial t}$ and
the previous inequality becomes 
$$
\displaystyle\frac{1}{2\alpha }\left| -\delta \otimes J\displaystyle\frac{%
\partial u}{\partial t}\right| ^{2}\leq \left( -\delta \otimes J\displaystyle%
\frac{\partial u}{\partial t},u\right) +\beta +\gamma .\eqno(6) 
$$
In the hypothesis' conditions, by integration of the inequality (6) we have 
\[
\displaystyle\frac{1}{2\alpha }\int_{T_{0}}\left| \displaystyle\frac{%
\partial u}{\partial t}\right| ^{2}dt^{1}\wedge ...\wedge dt^{p} 
\]
\[
+\int_{T_{0}}\left( \delta \otimes J\displaystyle\frac{\partial u}{\partial t%
},u\right) dt^{1}\wedge ...\wedge dt^{p}\leq \left( \beta +\gamma \right)
T^{1}...T^{p}. 
\]
By using the inequality from Theorem 3, we have 
\[
\displaystyle\frac{1}{2\alpha }\int_{T_{0}}\left| \displaystyle\frac{%
\partial u}{\partial t}\right| ^{2}dt^{1}\wedge ...\wedge dt^{p} 
\]
\[
-\displaystyle\frac{\sqrt{p}\max_{i}\left\{ T^{i}\right\} }{2\pi }%
\int_{T_{0}}\left| \displaystyle\frac{\partial u}{\partial t}\right|
^{2}dt^{1}\wedge ...\wedge dt^{p}\leq \left( \beta +\gamma \right)
T^{1}...T^{p}. 
\]
So 
\[
\left( \displaystyle\frac{1}{2\alpha }-\displaystyle\frac{\sqrt{p}%
\displaystyle\max_{i}\left\{ T^{i}\right\} }{2\pi }\right)
\int_{T_{0}}\left| \displaystyle\frac{\partial u}{\partial t}\right|
^{2}dt^{1}\wedge ...\wedge dt^{p} 
\]
\[
\leq \left( \beta +\gamma \right) T^{1}...T^{p} 
\]
and, as consequence 
\[
\int_{T_{0}}\left| \displaystyle\frac{\partial u}{\partial t}\right|
^{2}dt^{1}\wedge ...\wedge dt^{p}\leq \displaystyle\frac{2\pi \alpha \left(
\beta +\gamma \right) T^{1}...T^{p}}{\pi -\alpha \displaystyle%
\max_{i}\left\{ T^{i}\right\} \sqrt{p}}. 
\]
By integration, the inequality 
\[
\delta \left| u\right| -\beta \leq H\left( t,u\right) 
\]
produces 
\[
\delta \int_{T_{0}}\left| u\left( t\right) \right| dt^{1}\wedge ...\wedge
dt^{p}-\beta T^{1}...T^{p}\leq \int_{T_{0}}H\left( t,u\right) dt^{1}\wedge
...\wedge dt^{p}. 
\]
Because $H\left( t,u\right) $ is convex in $u$, 
\[
H\left( t,u\right) -H\left( t,0\right) \leq \left( \nabla H\left( t,u\left(
t\right) \right) ,u\left( t\right) \right) , 
\]
we obtain 
\[
\int_{T_{0}}H\left( t,u\left( t\right) \right) dt^{1}\wedge ...\wedge dt^{p} 
\]
\[
\leq \int_{T_{0}}\left[ H\left( t,0\right) +\left( \nabla H\left( t,u\left(
t\right) \right) ,u\left( t\right) \right) \right] dt^{1}\wedge ...\wedge
dt^{p} 
\]
\[
\leq \gamma T^{1}...T^{p}-\int_{T_{0}}\left( \delta \otimes J\frac{\partial u%
}{\partial t},u\left( t\right) \right) dt^{1}\wedge ...\wedge dt^{p} 
\]
\[
\leq \gamma T^{1}...T^{p}+\frac{\sqrt{p}\displaystyle\max_{i}\left\{
T^{i}\right\} }{2\pi }\int_{T_{0}}\left| \frac{\partial u}{\partial t}%
\right| ^{2}dt^{1}\wedge ...\wedge dt^{p} 
\]
\[
\leq \gamma T^{1}...T^{p}+\frac{\sqrt{p}\displaystyle\max_{i}\left\{
T^{i}\right\} 2\pi \alpha \left( \beta +\gamma \right) T^{1}...T^{p}}{2\pi
\left( \pi -\alpha \displaystyle\max_{i}\left\{ T^{i}\right\} \sqrt{p}%
\right) }. 
\]
By consequence 
\[
\int_{T_{0}}\left| u\left( t\right) \right| dt^{1}\wedge ...\wedge dt^{p} 
\]
\[
\leq \frac{1}{\delta }\left( \beta T^{1}...T^{p}+\gamma T^{1}...T^{p}+\frac{%
\sqrt{p}\displaystyle\max_{i}\left\{ T^{i}\right\} \alpha \left( \beta
+\gamma \right) }{\pi -\alpha \displaystyle\max_{i}\left\{ T^{i}\right\} 
\sqrt{p}}T^{1}...T^{p}\right) , 
\]
meaning that 
\[
\int_{T_{0}}\left| u\left( t\right) \right| dt^{1}\wedge ...\wedge
dt^{p}\leq \frac{\left( \beta +\gamma \right) T^{1}...T^{p}\pi }{\delta
\left( \pi -\alpha \displaystyle\max_{i}\left\{ T^{i}\right\} \sqrt{p}%
\right) } 
\]
and the proof ends.

\ 

\centerline{\bf References}

\ \ \ \ \ \ \ \ \ \ \ \ \ \ \ \ \ \ \ \ \ \ 

[1] I. Duca, A-M. Teleman, C. Udri\c{s}te: {\it Poisson-Gradient Dynamical
Systems with Convex Potential}, Proceedings of the 3-rd International
Colloquium '' Mathematics in Engineering and Numerical Physics '', 7-9
October, 2004, Bucharest.

[2] M. Forger, C. Paufler, H. Romer, {\it The Poisson bracket for Poisson
forms in multisymplectic field theory}, Reviews in Mathematical Physics, 15,
7 (2003), 705-743.

[3] I. V. Kanatchikov: {\it Geometric (pre)quantization in the
polysymplectic approach to field theory}, arXiv: hep-th/0112263 v3, 3 Jun
2002, 1-12; Differential Geometry and Its Applications, Proc. Conf., Opava
(Czech Republic), August 27-31, 2001, Silesian University Opava, 2002.

[4] L. V. Kantorovici, G. P. Akilov: {\it Analiz\u{a} func\c{t}ional\u{a}},
Editura \c{s}tiin\c{t}ific\u{a} \c{s}i enciclopedic\u{a}, Bucure\c{s}ti,
1980.

[5] J. Mawhin, M. Willem: {\it Critical Point Theory and Hamiltonian Systems}%
, Springer-Verlag, 1989.

[6] S. G. Mihlin: {\it Ecua\c{t}ii liniare cu derivate par\c{t}iale},
Editura \c{s}tiin\c{t}ific\u{a} \c{s}i enciclopedic\u{a}, Bucure\c{s}ti,
1983.

[7] M. Neagu, C. Udri\c{s}te : {\it From PDE \ Systems and Metrics to
Geometric Multi-time Field Theories}, Seminarul de Mecanic\u{a}, Sisteme
Dinamice Diferen\c{t}iale, 79, Universitatea de Vest din Timi\c{s}oara, 2001.

[8] P. J. Olver: \ {\it The\ \ Equivalence \ Problem and Canonical Forms for
Quadratic Lagrangians}, Advances in Applied Mathematics, 9, (1998), 226-227.

[9] C. Paufler, H. Romer: {\it De Donder-Weyl equations and multisymplectic
geometry}, arXiv: math-ph/0506022 v2, 20 Jun 2005.

[10] N. Roman-Roy, {\it Multisymplectic Lagrangian and Hamiltonian formalism
of first -order classical field theories}, arXiv: math-ph/0107019 v1, 20 Jul
2001, vol XX (XXXX), No. X, 1-9.

[11] R. E. Showalter: {\it Hilbert Space Methods for Partial Differential
Equations}, Electronic Journal of Differential Equations Monograph 01, 1994.

[12] A.-M. Teleman, C. Udri\c{s}te: {\it On Hamiltonian Formalisms in
Mathematical Physics}, 4-th International Conference of Balkan Society of
Geometers, Aristotle University of Thessaloniki, 26-30 June, 2002.

[13] C. Udri\c{s}te: {\it Nonclassical Lagrangian Dynamics and Potential Maps%
}, Conference in Mathematics in Honour of Professor Radu Ro\c{s}ca on the
occasion of his Ninetieth Birthday, Katholieke University Brussel,
Katholieke University Leuwen, Belgium, December 11-16, 1999; \newline
http://xxx.lanl.gov/math.DS/0007060, (2000).

[14] C. Udri\c{s}te: {\it Geometric Dynamics}, Kluwer Academic Publishers,
Dodrecht /Boston /London, Mathematics and Its Applications, 513, (2000);
Southeast Asian Bulletin of Mathematics, Springer-Verlag 24 (2000).

[15] C. Udri\c{s}te: {\it Solutions of DEs and PDEs as Potential Maps Using
First Order Lagrangians}, Centenial Vr\^{a}nceanu, Romanian Academy,
University of Bucharest, June 30-July 4, (2000);
http://xxx.lanl.gov/math.DS/ 0007061, (2000); Balkan Journal of Geometry and
Its Applications 6, 1, 93-108, 2001.

[16] C. Udri\c{s}te: {\it From integral manifolds and metrics to potential
maps}, Atti del'Academia Peloritana dei Pericolanti,Classe 1 di Scienze Fis.
Mat. e Nat., 81-82, A 01006 (2003-2004), 1-14.

[17] C. Udri\c{s}te, I. Duca: {\it Periodical Solutions of Multi-Time
Hamilton Equations,} Analele Universita\c{t}ii Bucure\c{s}ti, 55, 1 (2005),
179-188.

[18] C. Udri\c{s}te, I. Duca: {\it Poisson-gradient Dynamical Systems with
Bounded Non-linearity}, manuscript.

[19] C. Udri\c{s}te, M. Neagu: {\it From PDE Systems and Metrics to
Generalized Field Theories}, http://xxx.lanl.gov/abs/math.DG/0101207.

[20] C. Udri\c{s}te, M. Postolache: {\it Atlas of Magnetic Geometric Dynamics%
}, Geometry Balkan Press, Bucharest, 2001.

[21] C. Udri\c{s}te, A.-M. Teleman: {\it Hamilton Approaches of Fields Theory%
}, IJMMS, 57 (2004), 3045-3056; ICM Satelite Conference in Algebra and
Related Topics, University of Hong-Kong , 13-18.08.02.

[21] C. Udri\c{s}te, M. Ferrara, D. Opri\c s, Economic Geometric Dynamics,
Geometry Balkan Press, Bucharest, 2004.

\ \ \ \ \ \ \ \ \ \ \ \ \ \ \ \ \ \ \ \ \ \ \ \ \ \ \ \ \ \ \ \ \ \ \ \ \ \
\ \ \ \ \ \ \ \ \ \ \ \ \ \ \ \ \ \ \ 

University Politehnica of Bucharest

Departament of Mathematics

Splaiul Independen\c{t}ei 313

060042 Bucharest, Romania

udriste@mathem.pub.ro

\end{document}